\documentclass[10pt,reqno,a4paper]{amsart}

\usepackage[pdftex]{graphicx}
\usepackage{amsmath}
\usepackage{amssymb,amsthm,hyperref,caption}
\usepackage{color}
\definecolor{meinBlau}{rgb}{0,0.2,0.65} 
\definecolor{blau}{rgb}{0,0,0.75} 
\definecolor{rot}{rgb}{0.74,0,0} 
\hypersetup{colorlinks,linkcolor=blau,citecolor=blue,urlcolor=meinBlau}
\usepackage{times}
\usepackage{enumitem}

\allowdisplaybreaks

\newtheorem{coroll}{Corollary}

\theoremstyle{definition}

\newtheorem{remark}{Remark}

\def\E{{\mathbb {E}}}

\newcommand{\N}{\ensuremath{\mathbb{N}}}

\DeclareMathOperator{\tilt}{tilt}

\DeclareMathOperator{\law}{\overset{\mathcal{L}}{=}}
\DeclareMathOperator{\claw}{\overset{\mathcal{L}}{\rightarrow}}

\begin{document}

\author[M.~Kuba]{Markus Kuba}
\address{Markus Kuba\\
Department Applied Mathematics and Physics\\
University of Applied Sciences - Technikum Wien\\
H\"ochst\"adtplatz 5, 1200 Wien} %
\email{kuba@technikum-wien.at}

\author[A.~Panholzer]{Alois Panholzer}
\address{Alois Panholzer\\
Institut f{\"u}r Diskrete Mathematik und Geometrie\\
Technische Universit\"at Wien\\
Wiedner Hauptstr. 8-10/104\\
1040 Wien, Austria} \email{Alois.Panholzer@tuwien.ac.at}

\title[Limit law of one-sided tree destruction]{ A note on the limit law of one-sided tree destruction}


\begin{abstract}
This short note serves an addendum to the article "Destruction of very simple trees" by Fill, Kapur and Panholzer (2004). 
Therein, the limit law of one-sided tree destruction was determined by its moment sequence. 
We add an identification of the limit law, using recent results of Bertoin (2022),
in terms of the local time of a noise reinforced Bessel process.
\end{abstract}

\maketitle

\section{Addendum: Identification of the limit law }
In Theorem 5.1 of~\cite{FKP2006} the limit law of the cost of root isolation $Y_n$ in the class of very simple trees, a subclass
of simply generated trees, was obtained. The random variable $Y_n$ itself satisfies the distributional equation
\[
Y_n \law Y_{K_n}+t_n\quad n\ge 2, 
\]
for some random variable $K_n$~\cite{FKP2006}, with initial value $Y_1=t_1$ and toll function $t_n=n^{a}$, with $a\ge 0$. It was shown that for a parameter $\sigma>0$ and $a'=a+\frac12$
it holds
\[
\frac{Y_n}{\sigma n^{a'}} \claw Y_{a'},
\]
where $Y_{\alpha'}$ is uniquely determined by its moment sequence $m_s=\E(Y^s_{a'})$, for positive integers $s\ge 1$:
\[
m_s=\frac{s!}{2^{s/2}}\prod_{k=1}^{s}\frac{\Gamma(ka')}{\Gamma(ka'+\frac12)}.
\]

Recent results of Bertoin~\cite{Bertoin2022} on the noise reinforced Bessel process and its local time allows us to observe the following Corollary to Theorem 5.1 of~\cite{FKP2006}.
\begin{coroll}
The limit law $Y_{a'}$ has the same distribution as a moment-shifted (also called size biased or moment-tilted) $\hat{L}_t$ of a local time of a noise reinforced 
Bessel process, scaled by $\frac1{\sqrt{2}}$:
\[
Y_{a'} \law \frac1{\sqrt{2}}\cdot \tilt_1\big(\hat{L_t}/\kappa(p,t)\big),
\]
with factor $\kappa(p,t)$ given in~\eqref{scale}, dimension $d=1$ and reinforcement parameter $p=\frac12-\frac1{4a'}$.

\smallskip

Additionally, $Y_{a'} \law \frac1{\kappa(p,1)\sqrt{2}} \cdot \hat{I}$, where the scaled exponential functional $\hat{I}$ is given by
\[
\hat{I} = \int_0^{\infty}\exp(-\frac12 \hat{\xi}_t)dt, 
\]
with $\hat{\xi}_t$ denoting the subordinator with Laplace-Exponent $\hat\Phi$, determined by 
\[
\hat\Phi(r)=2^{-\frac12}\big(\frac1{4a'}\big)^{\frac12}\frac{\Gamma(\frac12)}{\frac12 B(\frac12,4a'r)}.
\]
\end{coroll}

\smallskip

\begin{remark}[Moments of Gamma type]
For $a'=\frac12$ one obtains a Rayleigh distribution, as noted before~\cite{FKP2006}. 
For $m\in\N$ and $m\cdot a'=\frac12$, i.e. $a'=\frac14$, we observe cancellations and obtain interesting moments of Gamma type~\cite{Janson}.
For example, for $a'=\frac14$ we get
\[
m_s=\frac{1}{2^{s/2}}\frac{\Gamma(\frac14)\Gamma(\frac12)\Gamma(s+1)}{\Gamma(\frac{s+1}4)\Gamma(\frac{s+2}4)},\quad s\ge 1.
\]
Density functions for such moment sequences can be obtained using inverse
Mellin computations; compare with ~\cite[Proof of Theorem 4.1]{BKW2021} or see~\cite[Theorem 5.4]{Janson}.
\end{remark}

\smallskip

\begin{proof}
The identification is directly done using moment sequences appearing in~\cite{Bertoin2022}. There, a noise reinforced Bessel process of dimension $d > 0$ and
with reinforcement parameter $p\in(-\infty,1/2)$ was considered. Set
\[
\alpha=1-d/2\in(0,1) \text{ and } \beta=\frac{\alpha}{1-2p}, 
\]
or equivalently
\[
p=\frac12 - \frac{\alpha}{2\beta},\text{ with }\beta >0.
\]
The local time $\hat{L}_t$ of the noise reinforced Bessel process has power moments (see \cite[Theorem 1.2]{Bertoin2022}), which can be written as
\[
\E(\hat{L_t}^s)=\kappa(p,t)^s\cdot \frac{1-2p}{\Gamma(1+\alpha)}\cdot \Gamma(s)\cdot\prod_{j=1}^{s-1}\frac{\Gamma(j\beta)}{\Gamma(\alpha+j\beta)}.
\]
with scale factor given by
\begin{equation}
\label{scale}
\kappa(p,t)=\frac{(2t)^{\alpha}\cdot\Gamma(1+\alpha)}{(1-2p)^{\alpha}\cdot\Gamma(1-\alpha)}.
\end{equation}
This implies that the scaled random variable $L=\hat{L_t}/\kappa(p,t)$ has moment sequence 
$\mu_s=\E(L^s)$, given by
\[
\mu_s=\frac{1-2p}{\Gamma(1+\alpha)}\cdot \Gamma(s)\cdot\prod_{j=1}^{s-1}\frac{\Gamma(j\beta)}{\Gamma(\alpha+j\beta)}, \quad s\ge 1.
\]
Using standard results on moment shifts (see for example~\cite[Lemma 3.4]{BKW2021}), also called size biased distributions or tilted distributions, 
we observe that the random variable $T=\tilt(L)$ has moment sequence $\E(T^s)=\frac{\mu_{s+1}}{\mu_1}$, 
or more explicitly 
\[
\E(T^s)=\Gamma(s+1)\cdot\prod_{j=1}^{s}\frac{\Gamma(j\beta)}{\Gamma(\alpha+j\beta)}.
\]
Since $\Gamma(s+1)=s!$, we finally set $\beta=a'$ and $\alpha=\frac12$. This leads to the desired moment sequence. 
Finally, we note that the choices lead to a dimension $d=1$ and $p=\frac12-\frac1{4a'}$, such that $1-2p=\frac1{4a'}$.

\smallskip 

Concerning the second statement we use the definition of $\hat{I}$ and its Laplace exponent~\cite[Equations 4.3 and 4.4]{Bertoin2022}. 
We use Corollary 4.3 of~\cite{Bertoin2022}, which relates $\hat{I}$ and $\hat{L}_1$. Note that there is a small misprint~\cite{Bertoin2022c} in the stated equation, as the prefactor
\[
c=(1/2-p)^{-\alpha}\frac{1-2p}{\Gamma(1-\alpha)},
\]
which can be rewritten into 
\[
\Big(\frac{1-2p}{2}\Big)^{-\alpha}\frac{1-2p}{\Gamma(1-\alpha)}
=\frac{2^{\alpha}(1-2p)^{1-\alpha}}{\Gamma(1-\alpha)}=\E(\hat{L}_1),
\]
should be replaced by its reciprocal. We state here the corrected version: 
\[
\E\big(f(\hat{I})\big)=\frac{1}{\E(\hat{L}_1)}\cdot \E\big(\hat{L}_1f(\hat{L}_1)\big).
\]
Choosing $f$ as the power functions, we observe a tilt of the moments~\cite[Lemma 3.4]{BKW2021}, 
leading to the moment sequence of $\hat{I}$ as the shifted moments of $\hat{L}_1$.
\end{proof}

\begin{remark}[More moments of Gamma type]
As already pointed out~\cite{Bertoin2022}, for $p=0$, such that $\alpha=\beta$, the random variable $\hat{L_t}$ 
simplifies to a Mittag-Leffler distribution:
\[
\E(\hat{L_t}^s)=\kappa(0,t)^s\cdot \frac{1}{\Gamma(1+\alpha)}\cdot \Gamma(s)\cdot\frac{\Gamma(\alpha)}{\Gamma(s\alpha)}.
=\kappa(0,t)^s\frac{\Gamma(s+1)}{\Gamma(s\alpha+1)}.
\]

\smallskip

The random variable $T$ also leads to moments of Gamma types for small $\beta$, proportional to $\alpha$:
$\beta=\alpha/m$, $m\in\N$, leading to moments
\[
\E(T^s)=\Gamma(s+1)\cdot\prod_{j=1}^{m}\frac{\Gamma(\frac{j\alpha}{m})}{\Gamma(\frac{(s+j)\alpha}{m})},\quad s\ge 1.
\]
Again, density functions can be obtained for such moments.
\end{remark}

\section{Outlook and Acknowledgments}
The authors are currently investigating into other occurrences of such limit laws
in combinatorial probability. 

\smallskip

The authors warmly thank Jean Bertoin for feedback on his work and explanations.

\end{document}